\documentclass[reprint,amsmath,amssymb,aps,pra,nofootinbib]{revtex4-2}

\usepackage[T1]{fontenc}
\usepackage[utf8]{inputenc}
\usepackage[english]{babel}
\usepackage[margin=1.8cm]{geometry}
\usepackage{graphicx}
\usepackage{ragged2e}
\usepackage{booktabs}
\usepackage{makecell}
\usepackage{fancyhdr}
\usepackage{titlesec}
\usepackage[unicode]{hyperref}
\usepackage{orcidlink}
\usepackage{todonotes}

\renewcommand{\leq}{\leqslant}

\newcommand{\deq}[1]{\mathrel{\smash[t]{\overset{#1}{=}}}}
\DeclareMathOperator{\diag}{diag}
\DeclareMathOperator{\Rad}{Rad}
\DeclareMathOperator{\rank}{rank}
\DeclareMathOperator{\GL}{GL}
\DeclareMathOperator{\RM}{RM}
\DeclareMathOperator{\St}{St}
\newcommand{\F}{\mathbb{F}}

\usepackage{titlesec}
\makeatletter
\renewcommand\thesection{\arabic{section}}
\makeatother
\titleformat{\section}[block]{\normalfont\normalsize\scshape\centering}{\thesection.}{0.75em}{}
\titlespacing*{\section}{0pt}{\baselineskip}{0.5\baselineskip}
\titlespacing*{\subsection}{0pt}{\baselineskip}{0.25\baselineskip}

\fancypagestyle{plain}{
  \fancyhf{}
  \fancyfoot[C]{\thepage}

}

\begin{document}

\setlength{\abovedisplayskip}{3pt}
\setlength{\abovedisplayshortskip}{3pt}
\setlength{\belowdisplayskip}{3pt}
\setlength{\belowdisplayshortskip}{3pt}
\setlength{\headheight}{13pt}
\setlength{\parskip}{0pt}
\emergencystretch=1em

\title{Classification of Boolean Cubic Forms in Ten Variables}

\author{
Kirill Khoruzhii\orcidlink{0000-0003-4689-3812}$^{1,*}$,
Patrick Gel\ss\orcidlink{0000-0002-3645-9513}$^{1}$,
Sebastian Pokutta\orcidlink{0000-0001-7365-3000}$^{1,2}$
}
\affiliation{
$^1$Zuse Institute Berlin, Berlin, Germany\\
$^2$Technische Universit\"at Berlin, Germany
}

\begin{abstract}
We classify Boolean cubic forms in ten variables up to $\GL(10,2)$-equivalence. The catalog contains all $3\,691\,560$ nonzero orbits. For every orbit we provide a representative with small monomial count, the stabilizer order, and the alternating rank together with an explicit decomposition. The classification is obtained by rank-stratified enumeration. We verify completeness by the Burnside orbit count and independently by the orbit--stabilizer identity. We also provide a fast, complete $\GL(10,2)$-invariant. By polarization, this gives the first complete classification of alternating trilinear forms in dimension $10$ over $\F_2$.
\end{abstract}

\maketitle
\thispagestyle{fancy}

\begingroup
\renewcommand\thefootnote{\fnsymbol{footnote}}
\footnotetext[1]{khoruzhii@zib.de}
\endgroup

\section{Introduction} 

Whenever a finite group acts on a discrete object, many natural questions reduce from the object itself to the set of group orbits. For Reed--Muller codes the relevant symmetry is the action of $\GL(m,2)$ on $\F_2^m$, which descends to each quotient layer
\begin{equation*}
  \RM^*(r,m)=\RM(r,m)/\RM(r-1,m).
\end{equation*}
Coset weight enumerators~\cite{sugita_1996,markov_2025}, covering-radius computations~\cite{brier_2003,dougherty_2021,gao_2023}, and the structure of cosets up to affine equivalence~\cite{gillot_2023} all depend on a function only through its orbit, so one representative and one stabilizer order encode all data for an entire orbit. Realizing this reduction in practice requires explicit lists of representatives, which motivates this paper.

We focus on the third layer $\RM^*(3,m)$, the space of Boolean cubic forms, and address the case $m=10$ in this paper. By polarization (Sec.~\ref{sec:bcf}), Boolean cubic forms are in bijection with alternating trilinear forms \cite{hora_2021} on $\F_2^m$, and apart from the coding-theoretic uses above, the same objects appear in the geometry of trivectors~\cite{aschbacher_1990,draisma_2010}, in the classification of code loops via their associator forms~\cite{obrien_2017}, in cryptographic equivalence problems for alternating tensors~\cite{beullens_2023}, and in the synthesis of fault-tolerant quantum circuits, where the alternating rank of the cubic part of a phase polynomial controls the Toffoli count~\cite{khoruzhii_2026}.

The classification of alternating trilinear forms has been carried out in several dimensions $m$ over various base fields. For $m \leq 8$ complete classifications are available over arbitrary fields~\cite{cohen_1988, noui_1997, midoune_2013,hora_2015}. For $m=9$, classifications are available over $\mathbb{C}$~\cite{vinberg_1978}, $\mathbb{R}$~\cite{borovoi_2022}, and $\F_2$~\cite{brier_2003,hora_2021}. To our knowledge, the case $m=10$ has not been addressed over any field.

We give the first complete classification in dimension $m=10$ over $\F_2$, in the form of an explicit catalog of all $3\,691\,560$ nonzero $\GL(10,2)$-orbits of Boolean cubic forms. For each orbit the catalog records the order of its stabilizer in $\GL(10,2)$ and a representative with small monomial count. We also construct a complete $\GL(10,2)$-invariant, encoded as a single $64$-bit word and verified to take distinct values on all representatives. The orbit count agrees with the Burnside summation~\cite{hou_1996}, and the recorded stabilizer orders provide a second independent check of completeness.

Alternating rank provides both the search coordinate and one of the useful outputs of the classification. It is the minimum number of decomposable cubics $u(x)v(x)w(x)$ whose sum is the given form. We use the resulting rank stratification to explore the orbit space: starting from the zero form, the enumeration proceeds by adding decomposable cubics and retaining the new $\GL(10,2)$-orbits reached at each layer. For every orbit we compute the alternating rank and store a decomposition realizing it.

The main orbit-identification tool is an invariant built from the finite-difference geometry of a form. We start from the orthogonality graph $G_f$ introduced in~\cite{hora_2021}. We verify that $G_f$ is already a complete invariant for $m<10$, but in dimension $10$ it has collisions. To resolve these, we introduce a bipartite incidence graph $B_f$ from the same polar bilinear forms, complementing the kernel information in $G_f$ with image information. We show that graph-isomorphism testing for the pair $(G_f,B_f)$ separates all representatives. The catalog uses a faster replacement: a fixed collection of local statistics extracted from these graphs, encoded as a complete $\GL(10,2)$-invariant in a single $64$-bit word.

\begin{figure*}[t]
    \centering
    \includegraphics[width=0.9\textwidth]{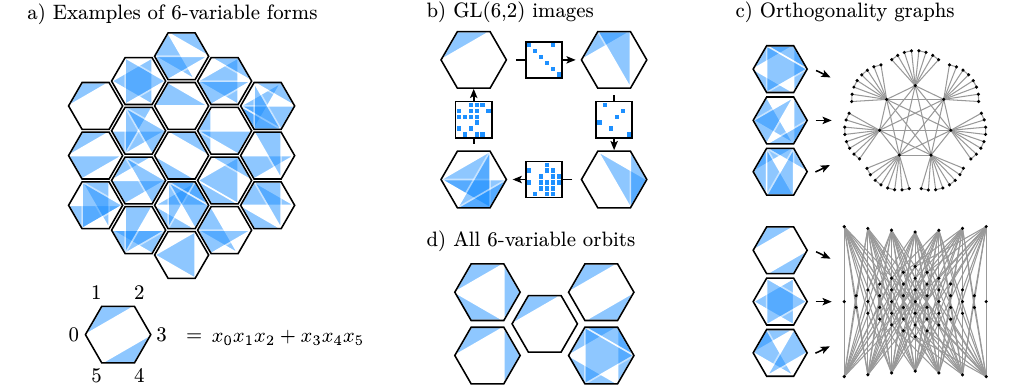}
    \caption{\justifying
        \textbf{Boolean cubic forms in six variables.}
        a) Examples of forms on $\F_2^6$, drawn as hypergraphs on a labeled hexagon; each triangular face corresponds to one monomial $x_i x_j x_k$.
        b) Examples of the $\GL(6,2)$-action. The basis change is displayed as a binary matrix, with blue squares for $1$ and empty squares for $0$.
        c) Orthogonality graphs $G_f$ for two orbits; equivalent forms determine isomorphic graphs.
        d) One representative of each of the five nonzero $\GL(6,2)$-orbits.
    }
    \label{fig:n6}
\end{figure*}

\section{Boolean Cubic Forms} \label{sec:bcf}

A Boolean cubic form on $\F_2^m$ is a polynomial~\cite{brier_2003, hora_2021}
\begin{equation*}
    f(x) = \sum_{ijk} T_{ijk}\, x_i x_j x_k,
\end{equation*}
with $0\leq i<j<k<m$ and $T_{ijk}\in \F_2$. We treat $f$ as a function $\F_2^m \to \F_2$, so all polynomials are reduced modulo $x_j^2 = x_j$ and written in square-free form. For polynomials $g, h$, we write
\begin{equation}
    g \deq{k} h
    \label{eq:deq}
\end{equation}
when their homogeneous parts of degree $k$ coincide. Boolean cubic forms make up the quotient $\RM^*(3,m)$, whose dimension is the binomial coefficient $\binom{m}{3}$.

For $u \in \F_2^m$, the \emph{finite difference}~\cite{brier_2003} of $f$ is
\begin{equation*}
    \Delta_u f(x) = f(x+u) + f(x),
\end{equation*}
a polynomial of degree one less than $f$. The third difference $T_f(u, v, w) = \Delta_u \Delta_v \Delta_w f$ is independent of $x$, trilinear and alternating in $(u, v, w)$. The map $f \mapsto T_f$ identifies $\RM^*(3,m)$ with the space of alternating trilinear forms on $\F_2^m$.

Forms that differ by a change of basis represent the same combinatorial object, so we classify forms up to the equivalence relation
\begin{equation*}
    f_1 \sim f_2 \ \Leftrightarrow \ \exists A \in \GL(m,2): f_1(x) \deq{3} f_2(Ax),
\end{equation*}
where, as in~\eqref{eq:deq}, we require equality only of the cubic part, ignoring the constant, linear, and quadratic terms. The set of nonzero equivalence classes is 
\begin{equation*}
    \mathcal{O}_m = \bigl(\RM^*(3,m)\setminus\{0\}\bigr)/\GL(m,2),
\end{equation*}
which is the main object of the classification. For example, $\mathcal{O}_6$ has five elements, shown in Fig.~\ref{fig:n6}d. Highly symmetric forms lie in small orbits. This is quantified by the \emph{stabilizer}
\begin{equation*}
    \St(f) = \{A \in \GL(m,2): f(Ax) \deq{3} f(x)\},
\end{equation*}
so the orbit of $f$ has size $|\GL(m,2)|/|\St(f)|$. The total count $|\mathcal{O}_m|$ follows from Burnside's lemma~\cite{hou_1996}; values for $m \leq 10$ are given in Tab.~\ref{tab:dim-summary}.

We use a compact notation for catalog representatives. A string $ijk$ denotes the monomial $x_i x_j x_k$, and sums of such strings are taken over $\F_2$; for example, $025+034$ means $x_0x_2x_5+x_0x_3x_4$. In the $m=10$ catalog, $f_j$ denotes the $j$-th representative in the catalog order.

A change of basis can eliminate one or more coordinates from $f$, leaving a form on a strict subspace of $\F_2^m$. These reducible directions $u$ form a subspace, the \emph{radical} of $f$,
\begin{equation*}
    \Rad(f) = \{u \in \F_2^m : \Delta_u f(x) \deq{2} 0\},
\end{equation*}
characterized by the vanishing of the quadratic part of the first difference~\cite{brier_2003}. The \emph{effective dimension}
\begin{equation*}
    \dim f = m - \dim\Rad(f)
\end{equation*}
is the number of coordinates $f$ genuinely uses, and $f$ is nondegenerate when $\Rad(f) = 0$.

For $u \in \F_2^m$, write $u(x) = \sum_j u_j x_j$. The \emph{alternating rank} (arank) of a cubic form $f$ is the smallest $r$ such that
\begin{equation*}
    f(x) \deq{3} \sum_{t=1}^r u_t(x) v_t(x) w_t(x),
\end{equation*}
with $u_t,v_t,w_t\in \F_2^m$. The alternating rank is invariant under the action of $\GL(m,2)$ and gives the rank stratification used below.

\begin{table*}[t]
\caption{Ambient-dimension summary. The orbit counts are for all nonzero Boolean cubic forms on $\F_2^m$.}
\label{tab:dim-summary}
\centering
\setlength{\tabcolsep}{8pt}
\begin{tabular}{lrrrrrrrr}
\toprule
$m$ & $3$ & $4$ & $5$ & $6$ & $7$ & $8$ & $9$ & $10$ \\
\midrule
$|\mathcal{O}_m|$ & $1$ & $1$ & $2$ & $5$ & $11$ & $31$ & $348$ & $3\,691\,560$ \\
$|\mathcal{R}_m|$ & $1$ & $15$ & $155$ & $1395$ & $11\,811$ & $97\,155$ & $788\,035$ & $6\,347\,715$ \\
$\dim \RM^*(3,m) = \binom{m}{3}$ & $1$ & $4$ & $10$ & $20$ & $35$ & $56$ & $84$ & $120$  \\
max alternating rank & $1$ & $1$ & $2$ & $3$ & $4$ & $5$ & $6$ & $7$ \\
max monomial count & $1$ & $1$ & $2$ & $5$ & $6$ & $8$ & $13$ & $17$ \\
\bottomrule
\end{tabular}
\end{table*}

\section{Enumeration by Alternating Rank} \label{sec:enumeration}

An alternating rank-$1$ cubic form $u(x)v(x)w(x)$ depends only on the three-dimensional subspace $\langle u,v,w\rangle\subseteq\F_2^m$, and vanishes when $u,v,w$ are linearly dependent. The set $\mathcal{R}_m$ of nonzero rank-$1$ forms is therefore indexed by three-dimensional subspaces of $\F_2^m$, and
\begin{equation*}
    |\mathcal{R}_m|=\binom{m}{3}_2=\frac{(2^m-1)(2^m-2)(2^m-4)}{(2^3-1)(2^3-2)(2^3-4)}.
\end{equation*}
The values for $m\leq 10$ are listed in Tab.~\ref{tab:dim-summary}. All elements of $\mathcal{R}_m$ are $\GL(m,2)$-equivalent: choosing $u,v,w$ as the first three basis vectors brings the form to $x_0x_1x_2$.

Adding elements of $\mathcal{R}_m$ defines a Cayley graph on $\RM^*(3,m)$, in which two forms are adjacent if they differ by one rank-$1$ form. The distance from $0$ in this graph is the alternating rank. Since $\GL(m,2)$ permutes $\mathcal{R}_m$, the graph is $\GL(m,2)$-equivariant, so this distance is well-defined on $\mathcal{O}_m$.

This yields the following BFS procedure. The first rank layer consists of a single orbit, represented by $x_0 x_1 x_2$. Each subsequent layer is generated by adding every rank-$1$ form to representatives of the previous layer; we identify the $\GL(m,2)$-orbit of each candidate via the invariants and keep only the new ones. Each orbit is expanded once, producing $|\mathcal{R}_m|$ candidates, so the total work is bounded by $|\mathcal{O}_m| \cdot |\mathcal{R}_m|$ invariant evaluations, far below the $2^{\binom{m}{3}}$ scale of the full space.

In practice the last expansion can be much smaller than the upper bound suggests. For $m=10$ we construct all layers up to rank $6$ exhaustively. To reach rank $7$, we sample $10^3$ rank-$6$ representatives and expand each of them. We use this sampled step both in the preliminary construction with the hyperplane-restriction profile $P_9$ described in App.~\ref{app:catalog-construction}, and in the final run with the complete invariant, which records alternating decompositions (Tab.~\ref{tab:rank-by-dim}). This sampling is used only to generate candidates: after orbit separation by the invariant below, reaching the Burnside count of~\cite{hou_1996} certifies that no orbit remains missing. The complete invariant is constructed in the next section.

\section{Complete Invariant} \label{sec:invariant}

The finite difference defines the \emph{orthogonality relation}, for $u,v\in\F_2^m$,
\begin{equation*}
    u \perp_f v \ \Leftrightarrow\ \Delta_u \Delta_v f(x) \deq{1} 0.
\end{equation*}
The orthogonality graph $G_f$ on $\F_2^m \setminus \{0\}$, introduced in~\cite{hora_2021}, has edges between distinct orthogonal pairs. The construction is $\GL(m,2)$-equivariant, so equivalent forms have isomorphic graphs (Fig.~\ref{fig:n6}c). For $m < 10$ the isomorphism class of $G_f$ is in fact a complete invariant of $f$, but at $m = 10$ this fails (see App.~\ref{app:marked-triples}).

To resolve the remaining collisions we introduce the incidence graph $B_f$, 
a bipartite graph with both parts indexed by $\F_2^m \setminus \{0\}$ and edges
\begin{equation*}
    u \to_f v \ \Leftrightarrow \ v(w) = 0 \text{ for all } w \perp_f u.
\end{equation*}
Left vertices act as vectors, right vertices as covectors. The construction 
is again $\GL(m,2)$-equivariant.

Both graphs come from the same object. The polar bilinear form $\beta_u (a,b) = \Delta_a \Delta_b \Delta_u f$ satisfies
\begin{equation*}
    \operatorname{ker} \beta_u = \{v : u \perp_f v\}, \ \ \ \  \operatorname{im} \beta_u = \{v : u \to_f v\}.
\end{equation*}
The neighborhoods of $u$ in $G_f$ and $B_f$ are thus the kernel and image of the same map.

Testing graph isomorphism directly is possible at this size, but it is not the most convenient way to build a fast invariant. Instead we use local graph statistics that are cheap to compute and sufficiently discriminative on the catalog. The two graphs above provide complementary statistics: closed-walk data on $G_f$, and right-side incidence data on $B_f$.

For the orthogonality graph $G_f$, let $A_f$ be its adjacency matrix. We use the closed-walk multisets
\begin{equation*}
    \mu_k(f)=\diag(A_f^k),
\end{equation*}
which record the number of closed walks of length $k$ based at each vertex. We verified on the complete catalogs for $m<10$ that $\mu_5$ separates all orbit representatives. In dimension $m=10$, it leaves $51$ unresolved collisions.

For the incidence graph, the useful extra information comes from the right side. For a right vertex $v$, let
\begin{equation*}
    N_L(v)=\{u:u\to_f v\}
\end{equation*}
be its left neighborhood. We define the right degree-codegree count $\theta_{s}^k(f)$ as the number of unordered pairs of distinct right vertices $v,v'$ such that $|N_L(v)|=|N_L(v')|=s$ and $|N_L(v)\cap N_L(v')|=k$.

Direct isomorphism testing of the incidence graphs $B_f$ separates all $51$ collisions. For a faster invariant, we use only a small fixed list of right degree-codegree counts. After the full $\RM^*(3,10)$ catalog was constructed, we treated the candidate counts $\theta_s^k$ as tests on the collision classes left by $\mu_5$ and selected a small set cover:
\begin{equation*}
    \theta = 
    \bigl(
        \theta_{216}^{56},\
        \theta_{222}^{56},\
        \theta_{228}^{64},\
        \theta_{231}^{40},\
        \theta_{237}^{68}
    \bigr).
\end{equation*}
Together with $\mu_5$, these counts separate all $m=10$ orbits except for the pair $f_{381},f_{382}$. The resolver $\rho$ separates this pair by the marked-triples test on $B_f$ (App.~\ref{app:marked-triples}).

By construction, the triple $(\mu_5,\theta,\rho)$ is a $\GL(10,2)$-invariant, so the number of distinct values it takes is at most the number of orbits. During enumeration we keep one nonzero representative per distinct value; reaching the Burnside count~\cite{hou_1996} shows this bound is attained with equality, so the triple in fact separates all $\GL(10,2)$-orbits. We compress it to a fixed $64$-bit hash word, verified collision-free on the resulting representatives, so two Boolean cubic forms in $10$ variables have the same word if and only if they are $\GL(10,2)$-equivalent.

\begin{table*}[t]
\caption{Number of nonzero orbits in $m = 10$, grouped by effective dimension $\dim f$ and alternating rank.}
\label{tab:rank-by-dim}
\centering
\setlength{\tabcolsep}{8pt}
\begin{tabular}{lrrrrrrrr}
\toprule
$\dim f$ & orbits & $\rank=1$ & $2$ & $3$ & $4$ & $5$ & $6$ & $7$ \\
\midrule
$3$ & $1$ & $1$ & $0$ & $0$ & $0$ & $0$ & $0$ & $0$ \\
$5$ & $1$ & $0$ & $1$ & $0$ & $0$ & $0$ & $0$ & $0$ \\
$6$ & $3$ & $0$ & $1$ & $2$ & $0$ & $0$ & $0$ & $0$ \\
$7$ & $6$ & $0$ & $0$ & $3$ & $3$ & $0$ & $0$ & $0$ \\
$8$ & $20$ & $0$ & $0$ & $2$ & $14$ & $4$ & $0$ & $0$ \\
$9$ & $317$ & $0$ & $0$ & $1$ & $21$ & $248$ & $47$ & $0$ \\
$10$ & $3\,691\,212$ & $0$ & $0$ & $0$ & $10$ & $1160$ & $2\,947\,440$ & $742\,602$ \\
\bottomrule
\end{tabular}
\end{table*}

\section{Results} \label{sec:results}

Applying the rank-stratified enumeration of Sec.~\ref{sec:enumeration}, with orbit identification by the invariant of Sec.~\ref{sec:invariant}, gives a catalog of all $3\,691\,560$ nonzero $\GL(10,2)$-orbits of Boolean cubic forms. In total $348$ orbits have nontrivial radical and are inherited from smaller effective dimensions, while the remaining $3\,691\,212$ orbits are nondegenerate in dimension $10$. The alternating rank distribution is shown in Tab.~\ref{tab:rank-by-dim}. The maximal alternating rank is $7$. Among nondegenerate $10$-dimensional forms only ranks $4,5,6,7$ occur, with rank $6$ forming the dominant layer.

We then applied the monomial-count minimization procedure of App.~\ref{app:thickness} to the representatives. The resulting forms are not certified to have minimum monomial count in their orbits, but they give explicit upper bounds and provide more compact representatives for display and reuse. The largest achieved count is $17$, as recorded in Tab.~\ref{tab:dim-summary}. Only one orbit remains at this value.

Stabilizer orders were computed for all catalog representatives by orbit--stabilizer backtracking, following the strategy of~\cite{hora_2021}. The distribution for the orbits is given in Tab.~\ref{tab:stab-dist}. There are $143$ distinct stabilizer orders, and most forms in $m=10$ have trivial stabilizer. The orbit--stabilizer identity is
\begin{equation*}
    \sum_f \frac{|\GL(10,2)|}{|\St(f)|}=2^{\binom{10}{3}}-1,
\end{equation*}
where the sum runs over the stored nonzero representatives. This gives a check, independent of the Burnside count~\cite{hou_1996}, that the represented orbits account for all nonzero Boolean cubic forms in ten variables.

All timings in this paragraph are measured on an AMD Ryzen 7 260. The complete $64$-bit invariant used during enumeration evaluates in $0.20$ ms per form on average. For comparison, we also tested direct graph-isomorphism certification by canonical labeling: once a canonical vertex order is fixed, the resulting colored adjacency matrix can be stored as a canonical hash. Direct canonical labeling with nauty~\cite{mckay_2014} did not finish within one minute for some of the orthogonality graphs $G_f$. We therefore used the augmented procedure of App.~\ref{app:gi}, which combines color-refinement signatures with component and block--cut decomposition. With this procedure, the canonical hash of $G_f$ takes $15$ ms per form on average, with worst case below $1.0$ s. For the colored incidence graph $B_f$, where the colors distinguish the two bipartite sides, the corresponding average is $5.5$ s per form; on the collision classes left by $\mu_5$, the worst case is $8.9$ s. Thus graph isomorphism is practical as a certification tool, but the enumeration relies on the faster local-statistics invariant.

The closed-walk multiset $\mu_5$ extracted from $G_f$ separates all representatives for $m<10$, hence the orthogonality graph $G_f$ itself is a complete invariant in these dimensions. In dimension $m=10$ this no longer holds: the pair $f_{381},f_{382}$ gives an explicit collision, with isomorphic orthogonality graphs but distinct $\GL(10,2)$-orbits. The additional statistics extracted from the incidence graph $B_f$ separate all collisions left by $\mu_5$. Thus the full colored pair $(G_f,B_f)$ is a complete $\GL(10,2)$-invariant.

The catalog data and the core code for computing the invariant and stabilizer orders are archived on Zenodo \href{https://zenodo.org/records/20773273}{\texttt{zenodo.org/records/20773273}}~\cite{bcf10_zenodo}. The accompanying repository, \href{https://github.com/khoruzhii/bcf10}{\texttt{github.com/khoruzhii/bcf10}}, contains supplementary code for the rank-stratified enumeration, graph-isomorphism certification, and auxiliary experiments.

\section{Discussion} \label{sec:discussion}

We expect the catalog to be useful for future Reed--Muller computations. Earlier classifications of cubic forms were used in this way for weight-enumerator and covering-radius questions in the third-order Reed--Muller codes~\cite{brier_2003,sugita_1996}. The present catalog provides the analogous input for $\RM^*(3,10)$: orbit representatives and stabilizer orders for all cubic cosets. Through the complementary map \cite{gillot_2023}, the same data also support computations in $\RM^*(7,10)$. The catalog further provides the orbit-indexing layer needed in Sarwate-type recursions~\cite{sarwate_1973,markov_2025} for the weight enumerator of $\RM(3,11)$.

In dimension $m=10$, the pair $(G_f,B_f)$ is complete on the classification, but it would be interesting to test whether $B_f$ alone already separates the $m=10$ orbits. For $m=11$, Burnside's formula gives about $6\cdot 10^{13}$ orbits~\cite{hou_1996}, so a catalog of the present kind would require substantially cheaper invariants and more structural use of restriction or projection profiles.

The monomial counts raise a separate extremal question. The minimum monomial count within a $\GL(m,2)$-orbit is one way to measure the algebraic complexity of a form. Our searches reduce every $m=10$ orbit to at most $17$ monomials, with only one orbit remaining at this value. Whether this candidate extremal orbit really has minimum monomial count equal to $17$, and hence whether the maximal value in $m=10$ is $17$, remains open.

The alternating rank decompositions point toward optimization applications. Although tensor-rank minimization is NP-hard in general~\cite{hillar_2013a}, the classification turns the ten-variable alternating case into a finite lookup problem: each orbit is represented together with an optimal decomposition. This gives a benchmark for rank-reduction heuristics such as flip graph search~\cite{kauers_2023,khoruzhii_2025,khoruzhii_2026}. It also suggests a local compilation strategy for quantum circuits: cubic phase-polynomial fragments on at most ten variables could be matched to catalog entries, up to a linear change of variables, and replaced by optimal decompositions. This is analogous in spirit to subcircuit replacement methods used for AND-count minimization~\cite{testa_2020a}.

\section*{Acknowledgments}

This research was supported by the DFG Cluster of Excellence MATH+ (EXC-2046/2, project id 390685689) funded by the Deutsche Forschungsgemeinschaft (DFG), as well as by the National High-Performance Computing (NHR) network.

\begin{table*}
\caption{Distribution of stabilizer orders among the $3\,691\,212$ nondegenerate $\GL(10,2)$-orbits of Boolean cubic forms.}
\label{tab:stab-dist}
\centering
\scriptsize
\setlength{\tabcolsep}{3pt}
\renewcommand{\arraystretch}{0.87}
\begin{minipage}[t]{0.495\textwidth}
\vspace{0pt}
\centering
\begin{tabular}{rrl}
\toprule
$|\St(f)|$ & \# orbits & representative \\
\midrule
21139292160 & 1 & 029+035+078+125+248+568 \\
18119393280 & 1 & 017+068+123+149+156 \\
15606743040 & 1 & 028+129+234+567 \\
1509949440 & 1 & 068+156+246+257+358+459+579+679 \\
1504051200 & 1 & \makecell[tl]{047+068+126+127+169+236+379+456+568 \\ 578} \\
402653184 & 1 & 026+038+047+059+137+235 \\
377487360 & 1 & 038+137+235+267+349+689 \\
301989888 & 2 & 038+158+246+369+678 \\
265420800 & 1 & 018+278+359+456 \\
201326592 & 1 & 028+035+069+168+346+378+589 \\
115605504 & 1 & 016+125+137+356+489 \\
61931520 & 1 & 034+125+128+159+256+578+689 \\
50331648 & 1 & 038+059+067+129+479+689 \\
37748736 & 1 & 024+038+059+067+125+278+356 \\
33554432 & 1 & 013+079+146+269+356+389 \\
25165824 & 2 & 029+057+135+189+459+568+679 \\
22118400 & 1 & 035+128+146+168+234+368+379 \\
18874368 & 2 & 059+078+138+149+268 \\
16777216 & 1 & 019+034+067+279+357+369+589 \\
13271040 & 1 & 014+089+124+236+257 \\
12386304 & 1 & 027+149+356+678 \\
9437184 & 2 & 016+089+248+345+378 \\
8847360 & 1 & 017+023+257+349+368 \\
8388608 & 1 & 014+179+239+245+346+357+478 \\
7077888 & 1 & 019+047+156+258+359 \\
6291456 & 2 & 018+149+157+256+346+458 \\
5898240 & 1 & 039+046+136+179+259+267+368+489+568 \\
5079040 & 2 & \makecell[tl]{025+029+034+038+067+069+089+124+168 \\ 236+289+456+478} \\
4194304 & 1 & 018+024+035+079+147+257+456+589 \\
3145728 & 10 & 012+035+048+268+348+379 \\
2752512 & 1 & 019+037+046+123+145+348+679 \\
2359296 & 1 & 024+089+139+156+239+479 \\
2088960 & 1 & 015+047+178+237+256+267+358+459+679 \\
2032128 & 1 & 014+039+137+258+346+679 \\
1572864 & 4 & 047+134+168+267+358+469 \\
1548288 & 1 & 056+079+134+258+267 \\
1474560 & 1 & 013+028+049+136+169+267+346+568 \\
1376256 & 2 & 059+127+146+256+257+349+789 \\
1179648 & 1 & 014+028+037+129+369+579 \\
1105920 & 1 & 049+058+156+245+249+289+379 \\
1048576 & 10 & 078+138+159+267+457+569 \\
887040 & 1 & \makecell[tl]{023+028+036+079+124+136+159+169+258 \\ 345+589+678} \\
884736 & 1 & 029+056+126+346+349+789 \\
786432 & 4 & 047+059+138+257+489+568 \\
774144 & 1 & 014+023+126+129+289+349+457+468 \\
737280 & 1 & 016+023+025+057+259+379+468 \\
688128 & 3 & 024+089+126+237+358+569 \\
655360 & 1 & \makecell[tl]{068+079+129+146+158+267+378+469+479 \\ 569} \\
589824 & 4 & 012+134+238+257+469+478 \\
524288 & 7 & 026+036+149+279+345+356+468 \\
442368 & 2 & 014+169+257+348+358 \\
393216 & 5 & 017+023+245+279+356+589 \\
387072 & 1 & \makecell[tl]{037+046+056+058+129+136+148+157+238 \\ 247+248+256} \\
331776 & 1 & 039+158+248+249+267+489 \\
327680 & 1 & 025+067+128+348+469 \\
294912 & 6 & 067+129+268+356+478 \\
262144 & 15 & 027+035+128+246+569+678 \\
245760 & 1 & 049+057+127+156+167+258+379+678 \\
221184 & 1 & 019+025+048+147+168+179+349+789 \\
196608 & 11 & 023+047+146+159+178+458 \\
184320 & 3 & 019+027+038+068+126+169+345+679 \\
163840 & 1 & \makecell[tl]{026+058+128+168+236+247+357+368+569 \\ 789} \\
147456 & 7 & 013+019+237+389+459+678 \\
138240 & 2 & \makecell[tl]{014+025+035+048+059+068+124+136+159 \\ 237+249+268+379+789} \\
131072 & 15 & 016+038+059+246+258+567+789 \\
98304 & 8 & 023+079+158+267+268+348 \\
92160 & 1 & 019+035+128+146+236+238+348+579 \\
86016 & 1 & 029+036+069+146+237+258+389+569+678 \\
73728 & 9 & 023+148+169+257+289+346 \\
69120 & 1 & \makecell[tl]{027+028+034+038+126+135+145+189+247 \\ 279+379+389+457+468+489} \\
65536 & 22 & 013+028+049+126+157+345 \\
55296 & 1 & 025+047+078+148+237+269+358+469 \\
49152 & 17 & 012+256+369+458+479 \\
36864 & 9 & 025+035+126+258+348+479 \\
32768 & 17 & 038+045+178+237+369+479 \\
24576 & 20 & 014+089+129+236+378+568 \\
18432 & 10 & 023+145+169+246+789 \\
16384 & 37 & 029+038+134+268+459+679 \\
15360 & 1 & \makecell[tl]{013+038+048+056+128+145+167+178+269 \\ 349+358+789} \\
12288 & 23 & 026+079+148+259+278+346 \\
9216 & 11 & 038+146+237+256+359+789 \\
\bottomrule
\end{tabular}
\end{minipage}
\hfill
\renewcommand{\arraystretch}{0.888}
\begin{minipage}[t]{0.495\textwidth}
\vspace{0pt}
\centering
\begin{tabular}{rrl}
\toprule
$|\St(f)|$ & \# orbits & representative \\
\midrule
8192 & 56 & 059+078+127+169+346+379 \\
6144 & 42 & 017+039+169+248+347+356 \\
5376 & 1 & \makecell[tl]{023+049+056+126+139+145+156+247+346 \\ 357+379+589+679} \\
4608 & 12 & 028+057+169+247+345+378+489 \\
4096 & 58 & 035+126+149+379+457+468 \\
3456 & 2 & \makecell[tl]{016+023+047+068+135+167+236+259+347 \\ 459+469+578+789} \\
3072 & 68 & 018+147+237+259+356+479 \\
2688 & 1 & \makecell[tl]{013+015+028+057+146+179+237+356+389 \\ 459+678+689} \\
2304 & 8 & 027+138+169+268+269+345+468+589 \\
2048 & 142 & 012+047+189+347+369+458 \\
1920 & 1 & \makecell[tl]{012+059+138+139+157+158+245+246+269 \\ 279+347+369+458+579+678} \\
1536 & 104 & 018+049+129+267+345+479 \\
1344 & 5 & \makecell[tl]{014+027+037+126+135+145+289+367+369 \\ 459+478} \\
1152 & 22 & 015+038+169+247+357+458 \\
1024 & 161 & 012+038+179+234+239+356+458 \\
960 & 4 & \makecell[tl]{012+037+069+138+145+256+258+278+367 \\ 468+579+789} \\
864 & 3 & 012+157+347+369+468+589 \\
768 & 104 & 015+048+189+257+268+356+479 \\
720 & 3 & 028+034+069+145+167+237+259+389 \\
672 & 2 & \makecell[tl]{013+037+046+079+126+127+135+189+234 \\ 237+459+678} \\
640 & 1 & \makecell[tl]{012+026+039+057+134+138+156+189+236 \\ 248+358+459+679} \\
576 & 4 & 025+089+138+159+247+269+346+347+379 \\
512 & 274 & 034+127+149+235+369+568 \\
432 & 1 & \makecell[tl]{029+037+038+045+146+189+258+259+267 \\ 378+579} \\
384 & 115 & 016+159+239+258+347+679 \\
336 & 5 & 036+049+128+139+257+347+567+568 \\
320 & 1 & \makecell[tl]{013+016+038+057+126+149+239+248+256 \\ 367+458+678} \\
288 & 7 & 014+039+156+278+357+457+469 \\
256 & 434 & 013+057+136+239+267+459+468 \\
224 & 4 & \makecell[tl]{014+025+137+159+168+179+234+267+279 \\ 378+489+568} \\
192 & 122 & 017+026+089+128+239+358+469 \\
168 & 6 & 016+029+045+125+348+367+789 \\
160 & 1 & \makecell[tl]{027+039+046+128+159+169+247+347+389 \\ 578} \\
144 & 10 & 015+017+058+129+278+346+457+479+569 \\
128 & 474 & 016+078+238+256+345+479 \\
120 & 3 & \makecell[tl]{015+034+048+068+128+269+346+368+379 \\ 457} \\
96 & 98 & 026+038+057+124+237+569+789 \\
80 & 2 & \makecell[tl]{016+027+129+157+234+289+346+348+358 \\ 379+467+568} \\
72 & 6 & 019+057+136+145+168+236+278+349 \\
64 & 944 & 012+016+137+238+459+468+567 \\
60 & 7 & 012+036+047+126+148+237+269+567+589 \\
48 & 157 & 016+057+178+237+259+348+467 \\
36 & 7 & 026+089+134+179+237+458+569 \\
32 & 2367 & 017+025+049+136+248+567+789 \\
24 & 227 & 017+045+068+135+249+256+348+479 \\
21 & 7 & \makecell[tl]{016+024+059+128+134+269+357+389+457 \\ 458+479+678} \\
20 & 8 & \makecell[tl]{018+045+058+138+146+236+237+279+469 \\ 579} \\
18 & 17 & \makecell[tl]{025+034+078+124+159+167+268+356+379 \\ 489} \\
16 & 3086 & 015+048+056+134+269+378+579 \\
14 & 8 & \makecell[tl]{014+028+035+127+136+137+257+367+468 \\ 479+569} \\
12 & 324 & 015+049+129+148+267+356+378+468 \\
11 & 1 & \makecell[tl]{014+017+018+025+039+123+126+137+249 \\ 368+389+456+478+579} \\
10 & 22 & \makecell[tl]{012+018+027+049+156+189+259+289+346 \\ 378} \\
9 & 7 & \makecell[tl]{013+018+028+059+126+136+157+237+245 \\ 349+368+459+467+789} \\
8 & 4845 & 034+078+159+247+268+357+469 \\
7 & 12 & \makecell[tl]{019+046+047+138+236+248+357+457+569 \\ 789} \\
6 & 512 & 069+078+134+189+236+257+458 \\
5 & 34 & \makecell[tl]{039+057+068+125+134+169+239+246+289 \\ 367+478} \\
4 & 16869 & 012+067+159+234+356+378+489 \\
3 & 1174 & 025+037+089+147+158+239+268+469 \\
2 & 72187 & 012+037+046+134+249+359+458+678 \\
1 & 3585671 & 013+027+068+125+179+234+467+589 \\
\bottomrule
\end{tabular}
\end{minipage}
\end{table*}

\section{Appendix} \label{sec:appendix}

\subsection{Catalog Construction and Residual Searches} \label{app:catalog-construction}

Here we record how the catalog was constructed before the final compact invariant of Sec.~\ref{sec:invariant} was fixed. The main routing tool in the initial search was a hyperplane-restriction profile, also used in \cite{hora_2021}. For a nonzero vector $u\in\F_2^{10}$, let $f|_u$ denote the restriction of $f$ to the hyperplane $u(x)=0$. Concretely, choose an index $j$ with $u_j=1$, solve for $x_j$, substitute into $f$, reduce modulo $x_i^2=x_i$, and keep the cubic part. This gives a Boolean cubic form in nine variables, well defined up to $\GL(9,2)$-equivalence. Since $\mu_5$ is complete in $m=9$, we define the multiset
\begin{equation*}
    P_9(f)=\big\{\mu_5(f|_u) \colon u \in \F_2^{10} \setminus \{0\}\big\}.
\end{equation*}
This profile is not part of the final invariant: computing all $1023$ hyperplane restrictions is much more expensive than the staged invariant of Sec.~\ref{sec:invariant}. During the initial search, however, $P_9$ was strong enough to support the rank-stratified enumeration described in Sec.~\ref{sec:enumeration}. 

After the exhaustive BFS through rank $5$ and the sampled rank-$6$ expansion, the Burnside count showed that three orbits were missing. The orbit--stabilizer identity gives a useful constraint: for representatives $f$ of all nonzero orbits
\begin{equation*}
    \sum_f \frac{|\GL(m,2)|}{|\St(f)|}
    =
    2^{\binom{m}{3}}-1.
\end{equation*}
Thus, once all but three orbits are known, the stabilizer orders of the missing orbits must account for the remaining mass. Let $a_1,a_2,a_3$ be these stabilizer orders. The uncovered part of the sum gave
\begin{equation*}
    \frac{1}{a_1}+\frac{1}{a_2}+\frac{1}{a_3}
    =
    \frac{5159}{3\,047\,424}.
\end{equation*}
The natural first assumption was that the missing orbits had stabilizer orders already present in the partial catalog. Under this restriction the equation has the unique solution
\begin{equation*}
    \{a_1,a_2,a_3\}=\{960,1536,5\,079\,040\}.
\end{equation*}
This reduced the residual search to three highly constrained targets.

The final candidates were found by local searches around known representatives with the same stabilizer orders. We inspected one- and two-triple deformations and screened them by $|\St(f)|$ and $P_9(f)$. Candidates matching an existing representative in $P_9$ and $|\St(f)|$ were compared by graph isomorphism of $B_f$. This produced the three missing forms $f_{382}$, $f_{1097}$, and $f_{1386}$, which have the same $P_9$ profiles as $f_{381}$, $f_{1144}$, and $f_{1387}$, respectively, but have nonisomorphic $B_f$.

After adding the three residual representatives, the list matched the Burnside count of $3\,691\,560$ nonzero orbits. The final complete invariant $(\mu_5,\theta,\rho)$ takes distinct values on all representatives, giving the classification certificate used in the main text.

\subsection{Marked-Triples Test} \label{app:marked-triples}

Here we describe the faster version of the final resolver $\rho$ used for the remaining pair in Sec.~\ref{sec:invariant}: 
\begin{align*}
    f_{381} &= 025+029+034+038+067+069+089 \\
        &+124+168+236+289+456+478, \\
    f_{382} &= 045+069+089+148+159+246+258 \\
        &+289+349+356+479+678+789.
\end{align*}
Direct graph isomorphism testing on $B_f$ would also separate the pair, but the following test makes the difference visible inside the incidence graph itself.

For a left vertex $u$, write
\begin{equation*}
    N_R(u)=\{v:u\to_f v\}
\end{equation*}
for its right neighborhood. For the remaining pair, the left vertices of degree $15$ in $B_f$ form a set
\begin{equation*}
    U=\{u:|N_R(u)|=15\}.
\end{equation*}
In both cases $|U|=31$, and $U$ is the set of nonzero points of a $5$-dimensional subspace of $\F_2^{10}$. Fix any $u\in U$; the construction below is independent of this choice.

Let $\mathcal T_u$ be the set of independent triples
$T=\{a,b,c\}\subset U$ such that $u\in\langle a,b,c\rangle$ and
\begin{equation*}
    |N_R(a)\cap N_R(b)\cap N_R(c)|=7.
\end{equation*}
We split $\mathcal T_u=\mathcal T_u^{\mathrm{in}}\sqcup \mathcal T_u^{\mathrm{out}}$, where $\mathcal T_u^{\mathrm{in}}$ consists of the triples containing $u$. For $T\in\mathcal T_u^{\mathrm{in}}$, let $r(T)$ be the number of triples $T'\in\mathcal T_u^{\mathrm{out}}$ with $|T\cap T'| > 0$.  The maximum of $r(T)$ separates the pair: it is $4$ for $f_{381}$ and $5$ for $f_{382}$. Thus we define $\rho$ by
\begin{equation*}
    \rho(f)=
    \begin{cases}
        0, &
            \text{if } (\mu_5,\theta)(f) \neq (\mu_5,\theta)(f_{381}), \\
        \max\limits_{T \in \mathcal{T}_u^{\mathrm{in}}} r(T), & \text{otherwise},
    \end{cases}
\end{equation*}
so $\rho$ resolves the last collision left by $(\mu_5,\theta)$.

\subsection{Monomial-Count Minimization} \label{app:thickness}

For many purposes it is useful to replace an arbitrary orbit representative by an equivalent form with fewer monomials. Each monomial $x_i x_j x_k$ is a rank-$1$ cubic, so a representative with $s$ monomials gives an explicit alternating rank decomposition of length $s$. The minimum possible monomial count over a $\GL(10,2)$-orbit is known as the thickness of the orbit. We do not attempt to certify this minimum. The representatives stored in the catalog are obtained by heuristic minimization, and their monomial counts should therefore be read as achieved upper bounds on thickness.

The search uses transvections $x_i \leftarrow x_i+x_j$, which generate $\GL(10,2)$. Under this substitution, only monomials containing $i$ and not containing $j$ can change: each term $x_ix_px_q$ with $p,q\notin\{i,j\}$ toggles the corresponding term $x_jx_px_q$. Hence the move is evaluated by comparing two restricted slices. If $r$ is the set of active pairs $\{p,q\}$ in the $i$-slice, with $p,q\notin\{i,j\}$, and $d$ is the corresponding set in the $j$-slice, then the change in monomial count is
\begin{equation*}
    \Delta=|r|-2|r\cap d|.
\end{equation*}
Thus all $10\cdot 9$ directed transvections can be evaluated using only bit operations and popcounts.

Starting from a representative $f$, we perform greedy descent on this transvection graph. At each step all directed transvections are inspected, and one with maximal decrease in monomial count is applied, with random tie-breaking. When no decreasing move exists, we allow up to $20$ neutral moves with $\Delta=0$, using a tabu memory of length $8$ on directed moves to avoid immediate cycling. These neutral moves allow the search to cross small plateaus that would otherwise trap a purely greedy descent.

The local descent is embedded in an iterated local search. The first descent starts from the catalog representative. Each later run starts from the endpoint of the previous run after a random kick consisting of $10$ random transvections. We used $4000$ descents per orbit in the main pass.

This pass substantially reduces the size of the displayed forms. In dimension $9$, the largest achieved monomial count becomes $13$, improving on the largest count $36$ in the representatives of~\cite{hora_2021}. In dimension $10$, the largest achieved monomial count is $17$. Only one orbit remains at this value, represented by
\begin{align*}
f_{32534} &=012+013+046+057+079+089\\
          &+123+127+159+168+235+249\\
          &+367+378+456+478+678,
\end{align*}
with $|\St(f_{32534})|=3$. Additional searches did not reduce this representative.

\subsection{Graph Isomorphism} \label{app:gi}

The graph-isomorphism checks in Sec.~\ref{sec:invariant} were carried out by canonical labeling. For every graph we compute a canonical ordering of its vertices; when an isomorphism is claimed, we also recover the corresponding vertex permutation and verify it directly on the adjacency matrices. Canonical hashes are used only as compact fingerprints of the resulting canonical forms, not as standalone evidence for isomorphism or nonisomorphism.

The implementation follows the standard nauty-style individualization--refinement scheme~\cite{mckay_2014}. Starting from a vertex partition, it repeatedly refines cells by neighbor counts into active splitter cells. If the partition is not discrete, a target cell is individualized and the search branches over its vertices. Leaves are compared lexicographically by the resulting adjacency matrix, while discovered automorphisms are used to prune later branches. Vertex colors are supported throughout the refinement and canonicalization; for the incidence graph $B_f$ we use them to distinguish the two bipartite sides.

For the graphs arising here, the plain individualization--refinement search is not sufficient. We therefore add several problem-independent preprocessing and decomposition steps. The root partition is strengthened by local vertex signatures: degree, triangle count, bounded shell profiles, and an iterative Weisfeiler--Leman neighborhood hash~\cite{weisfeiler_1968}. Disconnected graphs are canonicalized componentwise, and the canonical components are sorted before being assembled. Finally, for connected graphs with articulation points we use the block--cut tree~\cite{arvind_2025}. Each biconnected block is canonicalized with colors encoding its adjacent articulation subtrees, and messages are propagated along the block--cut tree from its center. This reduces large tree-like parts of the search to canonicalization of much smaller colored blocks.

The canonical-labeling procedure was used for the graph-isomorphism claims that are not decided by the closed-walk statistics alone. For $m<10$, the multiset $\mu_5$ extracted from $G_f$ already separates all stored representatives, and therefore certifies completeness of the orthogonality graph in these dimensions. In dimension $10$, full canonicalization of $G_f$ shows that the orthogonality graph is no longer complete: it leaves $19$ nontrivial isomorphism classes, containing $47$ catalog representatives in total. The pair $f_{381},f_{382}$ is one explicit collision, discussed in App.~\ref{app:marked-triples}. Canonicalization of the colored incidence graphs $B_f$ separates all of these collisions, thus canonicalization of the pairs $(G_f,B_f)$ separates all stored representatives in the $m=10$ catalog.

\bibliography{references.bib}

\end{document}